\documentclass[12 pt]{article}
\usepackage[pagebackref]{hyperref}
\hypersetup{
colorlinks = true,
linkcolor = blue,
citecolor = blue}
\usepackage[utf8]{inputenc}
\usepackage[T1]{fontenc}

\usepackage[]{amsfonts}
\usepackage{amssymb}
\usepackage{amsmath,amsthm}
\def\couleur(#1 #2 #3)
	{
\special{ps::[end]}}

\def\bx#1{\setbox1=\hbox{\kern3pt{#1}\kern3pt}			
 \dimen1=\ht1 \advance\dimen1 by 3pt \dimen2=\dp1 \advance\dimen2 by 3pt
 \setbox1=\hbox{\vrule height\dimen1 depth\dimen2\box1\vrule}%
 \setbox1=\vbox{\hrule\box1\hrule}%
 \advance\dimen1 by .4pt \ht1=\dimen1
 \advance\dimen2 by .4pt \dp1=\dimen2 \box1\relax}

\def\wbb#1{\kern#1em}
\def\vci{\vrule  width.02em height1.47ex depth-.0ex}		
\def\11{{\rm\wbb{.2}\vci\wbb{-.37}1}}

\def\underset#1#2{\mathrel{\mathop{\kern0pt #2}\limits_{#1}}}

\def\overset#1#2{\mathrel{\mathop{\kern0pt #2}\limits^{#1}}}

\newtheorem{thm}{Theorem}[section]
\newtheorem{lem}[thm]{Lemma}
\newtheorem{cor}[thm]{Corollary}
\newtheorem{defin}[thm]{Definition}

\begin{document}

\title{Carleson measures and $H^{p}$ interpolating sequences  in the polydisc.}

\author{Eric Amar}

\date{}
\maketitle
 \renewcommand{\abstractname}{Abstract}

\begin{abstract}
Let $S$ be a sequence of points in ${\mathbb{D}}^{n}.$ Suppose
 that $S$ is $H^{p}$ interpolating. Then we prove that the sequence
 $S$ is Carleson, provided that $p>2.$ We also give a sufficient
 condition, in terms of dual boundedness and Carleson measure,
 for $S$ to be an $H^{p}$ interpolating sequence.\par 

\end{abstract}
\ \par 
\ \par 
\ \par 
\ \par 

\tableofcontents
\ \par 

\section{Introduction and definitions.}
\quad Recall the definition of Hardy spaces in the polydisc.\ \par 

\begin{defin}
The {\bf Hardy space} $H^{p}({\mathbb{D}}^{n})$ is the set of
 holomorphic functions $f$ in ${\mathbb{D}}^{n}$ such that, with
 $e^{i\theta }:=e^{i\theta _{1}}{\times}\cdot \cdot \cdot {\times}e^{i\theta
 _{n}}$ and $d\theta :=d\theta _{1}\cdot \cdot \cdot d\theta
 _{n}$ the Lebesgue measure on ${\mathbb{T}}^{n}$:\par 
\quad \quad \quad $\displaystyle {\left\Vert{f}\right\Vert}_{H^{p}}^{p}:=\sup _{r<1}\int_{{\mathbb{T}}^{n}}{\left\vert{f(re^{i\theta
 })}\right\vert ^{p}d\theta }<\infty .$
\end{defin}
The \textbf{Hardy space} $H^{\infty }({\mathbb{D}}^{n})$ is the
 space of holomorphic and bounded functions in ${\mathbb{D}}^{n}$
 equipped  with the sup norm.\ \par 
\quad The space $H^{p}({\mathbb{D}}^{n})$ possesses a reproducing kernel
 for any $a\in {\mathbb{D}}^{n},\ k_{a}(z)$:\ \par 
\quad \quad \quad $\displaystyle k_{a}(z)=\frac{1}{(1-\bar a_{1}z_{1})}{\times}\cdot
 \cdot \cdot {\times}\frac{1}{(1-\bar a_{n}z_{n})}.$\ \par 
And we have $\forall f\in H^{p}({\mathbb{D}}^{n}),\ f(a):={\left\langle{f,k_{a}}\right\rangle},$
 where ${\left\langle{\cdot ,\cdot }\right\rangle}$ is the scalar
 product of the Hilbert space $H^{2}({\mathbb{D}}^{n}).$\ \par 
For $a\in {\mathbb{D}}^{n},$ set $((1-\left\vert{a}\right\vert
 ^{2})):=(1-\left\vert{a_{1}}\right\vert ^{2})\cdot \cdot \cdot
 (1-\left\vert{a_{n}}\right\vert ^{2}).$ The $H^{p}$ norm of
 $k_{a}$ is ${\left\Vert{k_{a}}\right\Vert}_{p}=((1-\left\vert{a}\right\vert
 ^{2}))^{-1/p'}$ hence the  normalized reproducing kernel in
 $H^{p}({\mathbb{D}}^{n})$ is\ \par 
\quad \quad \quad $\displaystyle k_{a,p}(z)=\frac{(1-\left\vert{a_{1}}\right\vert
 ^{2})^{1/p'}}{(1-\bar a_{1}z_{1})}{\times}\cdot \cdot \cdot
 {\times}\frac{(1-\left\vert{a_{n}}\right\vert ^{2})^{1/p'}}{(1-\bar
 a_{n}z_{n})}.$\ \par 
\ \par 
\quad Let $S$ be a sequence of points in ${\mathbb{D}}^{n}.$ \ \par 
Let $\ell ^{0}(S)$ be the set of sequences on $S$ and define
 the restriction operator $R_{p}:H^{p}({\mathbb{D}}^{n})\rightarrow
 \ell ^{0}(S)$ by:\ \par 
\quad \quad \quad $\displaystyle \forall f\in H^{p}({\mathbb{D}}^{n}),\ R_{p}f:=\lbrace
 ((1-\left\vert{a}\right\vert ^{2}))^{1/p}f(a)\rbrace _{a\in S}.$\ \par 

\begin{defin}
We say that $S$ is $H^{p}({\mathbb{D}}^{n})$ {\bf interpolating}
 if $R_{p}(H^{p}({\mathbb{D}}^{n}))\supseteq \ell ^{p}(S).$
\end{defin}
We can state, for $1\leq p<\infty $:\ \par 

\begin{defin}
We say that the sequence $S$ is {\bf Carleson} if the operator
 $R_{p}$ is bounded from $H^{p}({\mathbb{D}}^{n})$ to $\displaystyle
 \ell ^{p}(S),$ i.e.\par 
\quad \quad \quad $\displaystyle \exists C>0,\ \forall f\in \ H^{p}({\mathbb{D}}^{n}),\
 {\left\Vert{R_{p}f}\right\Vert}_{\ell ^{p}(S)}\leq C{\left\Vert{f}\right\Vert}_{H^{p}({\mathbb{D}}^{n})}.$
\end{defin}

\begin{defin}

 Let $\mu $ be a Borel measure in ${\mathbb{D}}^{n}.$ We shall
 say that $\mu $ is a {\bf Carleson measure} if $\exists p,\
 1\leq p<\infty ,$ such that:\par 
\quad \quad \quad $\exists C>0,\ \forall f\in H^{p}({\mathbb{D}}^{n}),\ f\in L^{p}(\mu
 )$ and ${\left\Vert{f}\right\Vert}_{L^{p}(\mu )}\leq C{\left\Vert{f}\right\Vert}_{H^{p}({\mathbb{D}}^{n})}.$
\end{defin}

      For $z\in {\mathbb{D}}^{n}$ define the rectangle $R_{z}:=\lbrace
 \zeta \in {\mathbb{T}}^{n}::\left\vert{\zeta _{j}-z_{j}/\left\vert{z_{j}}\right\vert
 }\right\vert \leq 1-\left\vert{z_{j}}\right\vert ,\ j=1,...,n\rbrace .$\ \par 
Now, for any open set $\Omega $ in ${\mathbb{T}}^{n},$ the (generalised)
 Carleson region is $\Gamma _{\Omega }:=\lbrace z\in {\mathbb{D}}^{n}::R_{z}\subset
 \Omega \rbrace .$\ \par 

\begin{defin}
The measure $\mu $ is a {\bf rectangular Carleson measure} in
 ${\mathbb{D}}^{n}$ if:\par 
\quad \quad \quad $\exists C>0::\forall z\in {\mathbb{D}}^{n},\ \left\vert{\mu
 }\right\vert (\Gamma _{R_{z}})\leq C\left\vert{R_{z}}\right\vert .$
\end{defin}
\quad The "natural" generalisation of the Carleson embedding Theorem
 from the disc to the polydisc would be that: $\mu $ is a Carleson
 measure iff it is a rectangular Carleson measure. But Carleson~\cite{Carleson74}
 gave a counter example to this and A. Chang~\cite{AChang79}
 gave the following characterisation.\ \par 

\begin{thm}
The measure $\mu $ is a Carleson measure in ${\mathbb{D}}^{n}$
 iff, for any open set $\Omega $:\par 
\quad \quad \quad $\displaystyle \exists C>0::\forall \Omega ,\ \left\vert{\mu
 }\right\vert (\Gamma _{\Omega })\leq C\left\vert{\Omega }\right\vert .$
\end{thm}
Because this characterisation does not depend on $p,$ if $\mu
 $ is $p$-Carleson for a $1\leq p<\infty $ then it is $q$-Carleson
 for any $q.$ This justifies the absence of $p$ in the definition
 of Carleson measure.\ \par 
\ \par 
\quad For $S$ a sequence of points in ${\mathbb{D}}^{n},$ define the measure:\ \par 
\quad \quad \quad $\displaystyle \chi _{S}:=\sum_{a\in S}{((1-\left\vert{a}\right\vert
 ^{2}))\delta _{a}}.$\ \par 
Then we can see easily that the sequence $S$ is Carleson iff
 the measure $\chi _{S}$ is a Carleson measure.\ \par 
\ \par 
\quad In one variable the interpolating sequences were characterized
 by L. Carleson~\cite{CarlInt58} for $H^{\infty }({\mathbb{D}})$
 and by H. Shapiro and A. Shieds~\cite{ShapShields61} for $H^{p}({\mathbb{D}})$
 by the same condition:\ \par 
\quad \quad \quad $\displaystyle \forall b\in S,\ \prod_{a\in S,\ a\neq b}{d_{G}(a,b)}\geq
 \delta >0,$\ \par 
where $d_{G}(a,b):=\left\vert{\frac{a-b}{1-\bar ab}}\right\vert
 $ is the Gleason distance.\ \par 
\quad In several variables for the unit ball $\Omega ={\mathbb{B}}$
 or for the unit polydisc $\Omega ={\mathbb{D}}^{n},$ this characterisation
 is still an open question, even for $H^{2}(\Omega ).$\ \par 
\quad Nevertheless we have already some necessary conditions.\ \par 

\begin{thm}
(~\cite{Varo72}) If the sequence $S$ is $H^{\infty }({\mathbb{D}}^{n})$
 interpolating then the measure $\chi _{S}$ is rectangular Carleson.
\end{thm}
We shall need\ \par 

\begin{defin}
We say that the $H^{p}(\Omega )$ interpolating sequence $S$ has
 the {\bf linear extension property}, LEP, if there is a bounded
 linear operator $E:\ell ^{p}(S)\rightarrow H^{p}(\Omega ^{n})$
 such that ${\left\Vert{E}\right\Vert}<\infty $ and for any $\lambda
 \in \ell ^{p}(S),\ E(\lambda )(z)$ interpolates the sequence
 $\lambda $ in $H^{p}(\Omega )$ on $S.$
\end{defin}

\begin{thm}
(~\cite{amarIntPoly}) \label{p2}If the sequence $S$ is $H^{\infty
 }({\mathbb{D}}^{n})$ interpolating then the measure $\chi _{S}$
 is Carleson and for any $p\geq 1$ the sequence $S$ is $H^{p}({\mathbb{D}}^{n})$
 interpolating with the LEP.
\end{thm}
\quad For the ball we have a better result by P. Thomas~\cite{Thomas87}.
 See also~\cite{AmarWirtBoule07}.\ \par 

\begin{thm}
\label{p1}If the sequence $S$ is $H^{1}({\mathbb{B}})$ interpolating
 then the measure $\chi _{S}$ is Carleson.
\end{thm}
\quad The first main result of this work is an analogous result for
 the polydisc:\ \par 

\begin{thm}
Let $p>2$ and suppose that the sequence $S$ of points in ${\mathbb{D}}^{n}$
 is $H^{p}({\mathbb{D}}^{n})$ interpolating. Then the sequence $S$ is Carleson.
\end{thm}
\quad We also have some sufficient conditions.\ \par 
\quad Let $\Omega $ be the ball or the polydisc.\ \par 

\begin{thm}
(~\cite{Bernd85} for the ball; ~\cite{BernChanLin87} for the
 polydisc) Let $S$ be a sequence of points in $\Omega .$ Let
 (a), (b) and (c) denote the following statements:\par 
(a) There is a constant $\delta >0$ such that\par 
\quad \quad \quad $\displaystyle \forall b\in S,\ \prod_{a\in S,\ a\neq b}{d_{G}(a,b)}\geq
 \delta >0.$\par 
(b) $S$ is an interpolating sequence for $H^{\infty }(\Omega ).$\par 
(c) The sequence $S$ is separated and is a Carleson sequence.\par 
Then (a) implies (b), (b) implies (c). However, the converse
 for each direction is false for $n\geq 2.$
\end{thm}
\quad Now we shall need the following important definition.\ \par 

\begin{defin}
Let $S$ be a sequence of points in $\Omega .$ We say that $S$
 is a {\bf dual bounded sequence} in $H^{p}(\Omega )$ if there
 is a sequence $\lbrace \rho _{a}\rbrace _{a\in S}\subset H^{p}(\Omega
 )$ such that, with $\displaystyle k_{b,p'}$ the normalised reproducing
 kernel in $H^{p'}(\Omega )$ for the point $b\in \Omega $:\par 
\quad \quad \quad $\displaystyle \exists \ C>0,\ \forall a,b\in S,\ {\left\langle{\rho
 _{a},k_{b,p'}}\right\rangle}=\delta _{a,b}$ and $\ {\left\Vert{\rho
 _{a}}\right\Vert}_{H^{p}(\Omega )}\leq C.$
\end{defin}
\quad Using this definition we have, with $\Omega $ the ball ${\mathbb{B}}$
 or the polydisc ${\mathbb{D}}^{n}$:\ \par 

\begin{thm}
(~\cite{AmarExtInt06}) Let $S$ be a sequence of points in $\Omega
 .$ Suppose $S$ is $H^{p}(\Omega )$ dual bounded with either
 $p=\infty $ or $p\leq 2.$ Moreover suppose that $S$ is Carleson.
 Then $S$ is $H^{q}(\Omega )$ interpolating with the LEP for any $1\leq q<p.$
\end{thm}
\quad For the ball we have a better result.\ \par 

\begin{thm}
(~\cite{intBall09}) \label{CI10}Let $S$ be a sequence of points
 in ${\mathbb{B}}.$ Suppose $S$ is $H^{p}({\mathbb{B}})$ dual
 bounded. Then $S$ is $H^{q}({\mathbb{B}})$ interpolating with
 the LEP for any $1\leq q<p.$
\end{thm}
\ \par 
\quad The second main new result here is the extension of Theorem~\ref{CI10}
 to the bidisc.\ \par 

\begin{thm}
Let $S$ be a dual bounded sequence in $H^{p}({\mathbb{D}}^{2}),$
 such that the associated measure $\chi _{S}$ is Carleson. Then
 $S$ is $H^{s}({\mathbb{D}}^{2})$ interpolating with the LEP,
 for any $s\in \lbrack 1,p\lbrack .$
\end{thm}

\section{First main result.}
\quad We have the easy lemma.\ \par 

\begin{lem}
\label{IC4} Let $S$ be a dual bounded sequence in $H^{p}({\mathbb{D}}^{n}),$
 then for any $1\leq q\leq p,\ S$ is dual bounded in $H^{q}({\mathbb{D}}^{n}).$
\end{lem}
\quad Proof.\ \par 
To see this, just take $\rho _{a}:=\gamma _{a}k_{a,r}$ where
 $\gamma _{a}$ is the dual sequence in $H^{p}({\mathbb{D}}^{n}).$ Then,\ \par 
\quad \quad \quad \quad $\displaystyle {\left\langle{\rho _{a},k_{b,q'}}\right\rangle}={\left\langle{\gamma
 _{a}k_{a,r},k_{b,q'}}\right\rangle}=\gamma _{a}(b)k_{a,r}(b){\left\Vert{k_{b}}\right\Vert}_{p'}^{-1}.$\
 \par 
But\ \par 
\quad \quad \quad $\displaystyle {\left\langle{\gamma _{a},k_{b,p'}}\right\rangle}=\gamma
 _{a}(b){\left\Vert{k_{b}}\right\Vert}_{p'}^{-1}=\delta _{ab},$\ \par 
using\ \par 
\quad \quad \quad $\displaystyle \delta _{a,b}={\left\langle{\gamma _{a},k_{b,p'}}\right\rangle}={\left\Vert{k_{b}}\right\Vert}_{p'}^{-1}\gamma
 _{a}(b).$\ \par 
Recall that ${\left\Vert{k_{a}}\right\Vert}_{p}=((1-\left\vert{a}\right\vert
 ^{2}))^{-1/p'}$ hence defining $\chi _{a}:={\left\Vert{k_{a}}\right\Vert}_{2}^{-2}=((1-\left\vert{a}\right\vert
 ^{2}))$ we get ${\left\Vert{k_{a}}\right\Vert}_{p'}^{-1}=((1-\left\vert{a}\right\vert
 ^{2}))^{1/p}=\chi _{a}^{1/p}.$So we get\ \par 
\quad \quad \quad $\displaystyle {\left\langle{\rho _{a},k_{b,q'}}\right\rangle}=\gamma
 _{a}(b)k_{a,r}(b)\chi _{b}^{1/q}=\delta _{ab}\chi _{a}^{\frac{1}{q}-\frac{1}{p}}k_{a,r}(a)=\delta
 _{a,b}\chi _{a}^{\frac{1}{q}-\frac{1}{p}+\frac{1}{r'}}k_{a}(a).$\ \par 
But $k_{a}(a)={\left\Vert{k_{a}}\right\Vert}_{2}^{2}=\chi _{a}^{-1}$
 so finally:\ \par 
\quad \quad \quad $\displaystyle {\left\langle{\rho _{a},k_{b,q'}}\right\rangle}=\delta
 _{a,b}\chi _{a}^{\frac{1}{q}-\frac{1}{p}+\frac{1}{r'}-1}=\delta _{ab}$\ \par 
provided that $\frac{1}{r}=\frac{1}{q}-\frac{1}{p}$ which is
 possible because $1\leq q\leq p.$\ \par 
\quad It remains to prove that $\exists C>0::\forall a\in S,\ {\left\Vert{\rho
 _{a}}\right\Vert}_{H^{q}({\mathbb{D}}^{n})}\leq C.$ But because
 $\gamma _{a}\in H^{p}({\mathbb{D}}^{n})$ and $k_{a,r}\in H^{r}({\mathbb{D}}^{n})$
 and $\frac{1}{q}=\frac{1}{r}+\frac{1}{p}$ we get\ \par 
\quad \quad \quad $\displaystyle {\left\Vert{\rho _{a}}\right\Vert}_{q}\leq {\left\Vert{\gamma
 _{a}}\right\Vert}_{p}{\left\Vert{k_{a,r}}\right\Vert}_{r}\leq C$\ \par 
using that\ \par 
\quad \quad \quad $\displaystyle \forall a\in S,\ {\left\Vert{\gamma _{a}}\right\Vert}_{p}\leq
 C,\ {\left\Vert{k_{a,r}}\right\Vert}_{r}\leq C'.$\ \par 
This finishes the proof of the lemma. $\hfill\blacksquare $\ \par 
\ \par 
\quad Now we are in position to prove our first main result.\ \par 

\begin{thm}
\label{iDC7}Let $p>2$ and suppose that the sequence $S$ of points
 in ${\mathbb{D}}^{n}$ is $H^{p}({\mathbb{D}}^{n})$ interpolating.
 Then the sequence $S$ is Carleson.
\end{thm}
\quad Proof.\ \par 
Suppose that the sequence $S\subset {\mathbb{D}}^{n}$ is $H^{p}({\mathbb{D}}^{n})$
 interpolating.\ \par 
\quad Because $p\geq 2$ we have that $S$ is dual bounded in $H^{2}({\mathbb{D}}^{n})$
 by Lemma~\ref{IC4}. Call $\lbrace \rho _{a}\rbrace _{a\in S}$
 the dual sequence to the normalised reproducing kernels $k_{a,2}.$
 We have $\exists C>0::\forall a\in S,\ {\left\Vert{\rho _{a}}\right\Vert}_{2}\leq
 C.$\ \par 
Now we set $S_{N}$ the truncation of $S$ to its $N$ first terms
 and set $E_{S}:=\mathrm{S}\mathrm{p}\mathrm{a}\mathrm{n}\lbrace
 k_{a},\ a\in S\rbrace $ and $E_{S}^{N}:=\mathrm{S}\mathrm{p}\mathrm{a}\mathrm{n}\lbrace
 k_{a},\ a\in S_{N}\rbrace .$ \ \par 
Let also $P_{N}$ be the orthogonal projection from $H^{2}({\mathbb{D}}^{n})$
 onto $E_{S}^{N}.$ We have ${\left\Vert{P_{N}\rho _{a}}\right\Vert}_{2}\leq
 {\left\Vert{\rho _{a}}\right\Vert}_{2}.$\ \par 
Moreover the sequence $\lbrace P_{N}\rho _{a}\rbrace $ is still
 a dual sequence  to $\lbrace k_{a,2}\rbrace _{a\in S_{N}}$ because\ \par 
\quad \quad \quad $\displaystyle \forall a,b\in S,\ {\left\langle{P_{N}\rho _{a},k_{b,2}}\right\rangle}={\left\langle{\rho
 _{a},P_{N}k_{b,2}}\right\rangle}={\left\langle{\rho _{a},k_{b,2}}\right\rangle}=\delta
 _{a,b}.$\ \par 
Abusing the notation, we still denote by $\rho _{a}$ the dual
 sequence to $\lbrace k_{a,2}\rbrace _{a\in S_{N}}$ in $E_{S}^{N},$
 i.e. noting $\rho _{a}$ instead of $P_{N}\rho _{a}.$\ \par 
Let $\displaystyle \frac{1}{2}=\frac{1}{p}+\frac{1}{r}$ and $f\in
 H^{p}({\mathbb{D}}^{n}),\ g\in H^{r}({\mathbb{D}}^{n}).$ Because
 $fg\in H^{2}({\mathbb{D}}^{n}),$ ${\left\Vert{P_{N}(fg)}\right\Vert}_{2}\leq
 {\left\Vert{fg}\right\Vert}_{2},$ we have that:\ \par 
\quad \quad \quad \begin{equation}  P_{N}(fg)=\sum_{a\in S_{N}}{{\left\langle{fg,k_{a,2}}\right\rangle}\rho
 _{a}}=\sum_{a\in S_{N}}{f(a)\chi _{a}^{1/p}g(a)\chi _{a}^{1/r}\rho
 _{a}}.\label{ICG3}\end{equation}\ \par 
Let $\lbrace \epsilon _{a}\rbrace _{a\in S}$ be a random Bernouilli
 i.i.d.. Consider the finite sum $\sum_{a\in S_{N}}{\epsilon
 _{a}\mu _{a}\rho _{a}}.$ We have, fixing $\zeta \in {\mathbb{T}}^{n},$
 with $X(\zeta ):=\sum_{a\in S_{N}}{\epsilon _{a}\mu _{a}\rho
 _{a}(\zeta )}$:\ \par 
\quad \quad \quad $\displaystyle \mathrm{V}\mathrm{a}\mathrm{r}(X(\zeta ))={\mathbb{E}}(\left\vert{\sum_{a\in
 S_{N}}{\epsilon _{a}\mu _{a}\rho _{a}(\zeta )}}\right\vert ^{2})=\sum_{a\in
 S_{N}}{\left\vert{\mu _{a}}\right\vert ^{2}\left\vert{\rho _{a}(\zeta
 )}\right\vert ^{2}}$\ \par 
because the $\epsilon _{a}$ are independent.\ \par 
\quad By the Theorem of Fubini-Tonelli we get\ \par 
\quad \quad \quad $\displaystyle {\mathbb{E}}{\left({{\left\Vert{\sum_{a\in S_{N}}{\epsilon
 _{a}\mu _{a}\rho _{a}}}\right\Vert}_{H^{2}}^{2}}\right)}=\int_{T^{n}}{{\mathbb{E}}(\left\vert{\sum_{a\in
 S_{N}}{\epsilon _{a}\mu _{a}\rho _{a}(\zeta )}}\right\vert ^{2})d\zeta
 }=$\ \par 
\quad \quad \quad \quad \quad \quad \quad $\displaystyle =\int_{{\mathbb{T}}^{n}}{\sum_{a\in S_{N}}{\left\vert{\mu
 _{a}}\right\vert ^{2}\left\vert{\rho _{a}(\zeta )}\right\vert
 ^{2}d\zeta }}=\sum_{a\in S_{N}}{\left\vert{\mu _{a}}\right\vert
 ^{2}{\left\Vert{\rho _{a}}\right\Vert}_{2}^{2}}.$\ \par 
We always have that ${\left\Vert{\rho _{a}}\right\Vert}_{2}\geq
 1$ because $1={\left\langle{\rho _{a},k_{a,2}}\right\rangle}\leq
 {\left\Vert{\rho _{a}}\right\Vert}_{2}{\left\Vert{k_{a,2}}\right\Vert}_{2}={\left\Vert{\rho
 _{a}}\right\Vert}_{2}$ and because $S$ hence $S_{N}$ is dual
 bounded in $H^{2},$ we have that ${\left\Vert{\rho _{a}}\right\Vert}_{2}\leq
 C.$ So we get\ \par 
\quad \quad \quad $\displaystyle \sum_{a\in S_{N}}{\left\vert{\mu _{a}}\right\vert
 ^{2}}\leq {\mathbb{E}}{\left({{\left\Vert{\sum_{a\in S_{N}}{\epsilon
 _{a}\mu _{a}\rho _{a}}}\right\Vert}_{H^{2}}^{2}}\right)}\leq
 C\sum_{a\in S_{N}}{\left\vert{\mu _{a}}\right\vert ^{2}}.$\ \par 
\quad We notice that the constants here are independent of $N.$\ \par 
Hence we can find a sequence $\lbrace \epsilon _{a}\rbrace _{a\in
 S_{N}}$ such that\ \par 
\quad \quad \quad $\displaystyle {\left\Vert{\sum_{a\in S_{N}}{\epsilon _{a}\mu
 _{a}\rho _{a}}}\right\Vert}_{2}^{2}\geq \sum_{a\in S_{N}}{\left\vert{\mu
 _{a}}\right\vert ^{2}}$\ \par 
which means, see~\cite[p. 18]{AmarThesis77} for another way to get it,\ \par 
\quad \quad \quad \begin{equation}  \sup _{\lbrace \epsilon _{a}\rbrace }{\left\Vert{\sum_{a\in
 S_{N}}{\epsilon _{a}\mu _{a}\rho _{a}}}\right\Vert}_{2}^{2}\geq
 \sum_{a\in S_{N}}{\left\vert{\mu _{a}}\right\vert ^{2}}\label{ICG1}\end{equation}\
 \par 
\ \par 
\quad Now take $\lambda \in \ell ^{p}(S)$ and for $\lbrace \epsilon
 _{a}\rbrace _{a\in S}$ a random Bernouilli i.i.d., we can choose
 a $f\in H^{p}({\mathbb{D}}^{n})$ such that ${\left\langle{f,k_{a,p}}\right\rangle}=\chi
 _{a}^{1/p}f(a)=\epsilon _{a}\lambda _{a}$ and ${\left\Vert{f}\right\Vert}_{H^{p}({\mathbb{D}}^{n})}\leq
 C{\left\Vert{\lambda }\right\Vert}_{\ell ^{p}(S)}$ because $S$
 is $H^{p}({\mathbb{D}}^{n})$ interpolating.\ \par 
\quad Applying~(\ref{ICG1}) to $P_{N}(fg)$ we get\ \par 
\quad \quad \quad $\displaystyle \mu _{a}:=f(a)\chi _{a}^{1/p}g(a)\chi _{a}^{1/r}=\epsilon
 _{a}\lambda _{a}g(a)\chi _{a}^{1/r},$\ \par 
and using~(\ref{ICG3}), \ \par 
\quad \quad \quad {\it  }$\displaystyle \sup _{\epsilon _{a}}{\left\Vert{\sum_{a\in
 S_{N}}{\epsilon _{a}\lambda _{a}g(a)\chi _{a}^{1/r}\rho _{a}}}\right\Vert}_{2}^{2}\geq
 \sum_{a\in S_{N}}{\left\vert{\lambda _{a}g(a)\chi _{a}^{1/r}}\right\vert
 ^{2}}.$\ \par 
So, setting $\nu _{a}:=g(a)\chi _{a}^{1/r},$ and $\nu :=\lbrace
 \nu _{a}\rbrace _{a\in S_{N}},$ we get by H\"older inequalities,
 using $\frac{1}{2}=\frac{1}{p}+\frac{1}{r},$\ \par 
\quad \quad \quad $\displaystyle {\left\Vert{\lambda }\right\Vert}_{\ell ^{p}}^{2}{\left\Vert{\nu
 }\right\Vert}_{\ell ^{r}}^{2}\geq \sum_{a\in S_{N}}{\left\vert{\nu
 _{a}\lambda _{a}}\right\vert ^{2}}$\ \par 
and\ \par 
\quad \quad \quad \begin{equation}  \sup _{\lambda \in \ell ^{p},{\left\Vert{\lambda
 }\right\Vert}_{\ell ^{p}}\leq 1}\sum_{a\in S_{N}}{\left\vert{\nu
 _{a}\lambda _{a}}\right\vert ^{2}}={\left\Vert{\nu }\right\Vert}_{\ell
 ^{r}}^{2}\label{ICG0}\end{equation}\ \par 
\quad So we choose $f::{\left\langle{f,k_{a,p}}\right\rangle}=\chi
 _{a}^{1/p}f(a)=\epsilon _{a}\lambda _{a}$ with $\lambda $ making
  ~(\ref{ICG0}) and $\lbrace \epsilon _{a}\rbrace _{a\in S_{N}}$
 making ~(\ref{ICG1}) (recall that $S_{N}$ is finite).  We get
 for this choice, with ${\left\Vert{f}\right\Vert}_{p}\leq C,$
 because ${\left\Vert{\lambda }\right\Vert}_{\ell ^{p}}\leq 1,$\ \par 
\quad \quad \quad $\displaystyle {\left\Vert{P_{N}(fg)}\right\Vert}_{2}^{2}={\left\Vert{\sum_{a\in
 S_{N}}{\epsilon _{a}\lambda _{a}g(a)\chi _{a}^{1/r}\rho _{a}}}\right\Vert}_{2}^{2}\geq
 \sum_{a\in S_{N}}{\left\vert{\lambda _{a}g(a)\chi _{a}^{1/r}}\right\vert
 ^{2}}.$\ \par 
Hence ${\left\Vert{\nu }\right\Vert}_{\ell ^{r}(S_{N})}\leq {\left\Vert{P_{N}(fg)}\right\Vert}_{2}\leq
 {\left\Vert{fg}\right\Vert}_{2}.$\ \par 
But ${\left\Vert{fg}\right\Vert}_{2}\leq {\left\Vert{f}\right\Vert}_{p}{\left\Vert{g}\right\Vert}_{r}$
 using $\frac{1}{2}=\frac{1}{p}+\frac{1}{r},$ hence\ \par 
\quad \quad \quad $\displaystyle {\left\Vert{\nu }\right\Vert}_{\ell ^{r(S_{N})}}\leq
 {\left\Vert{fg}\right\Vert}_{2}\leq {\left\Vert{f}\right\Vert}_{p}{\left\Vert{g}\right\Vert}_{r}\leq
 C{\left\Vert{g}\right\Vert}_{r}.$\ \par 
Now ${\left\Vert{\nu }\right\Vert}_{\ell ^{r}(S_{N})}^{r}=\sum_{a\in
 S_{N}}{\left\vert{g(a)\chi _{a}^{1/r}}\right\vert ^{r}}\leq
 C{\left\Vert{g}\right\Vert}_{r}^{r}$ so we get\ \par 
\quad \quad \quad $\displaystyle \sum_{a\in S_{N}}{\left\vert{g(a)}\right\vert
 ^{r}\chi _{a}}\leq C{\left\Vert{g}\right\Vert}_{r}^{r}.$\ \par 
Because none of the constants depends on $N,$ we get:\ \par 
\quad \quad \quad $\displaystyle \sum_{a\in S}{\left\vert{g_{a}}\right\vert ^{r}\chi
 _{a}}\leq C{\left\Vert{g}\right\Vert}_{r}^{r},$\ \par 
which means that $S$ is a Carleson sequence. $\hfill\blacksquare $\ \par 

\section{Second main result.}

\subsection{Interpolation of $H^{p}$ spaces.}
\quad We set the assumption:\ \par 
\quad \quad \quad (AS)     We have $(L^{1}({\mathbb{T}}^{n}),\ BMO({\mathbb{T}}^{n}))_{\theta
 }=L^{p}({\mathbb{T}}^{n}),$ where $\frac{1}{p}=1-\theta .$\ \par 
\ \par 
\quad We shall use the following theorem of Lin~\cite{Lin86} (see also~\cite{ChangFeff85},
 p. 37) proved for the bidisc:\ \par 

\begin{thm}
\label{AI10}If $1\leq p_{0}<p_{1}\leq \infty ,$ then $(H^{p_{0}}({\mathbb{T}}^{2}),\
 L^{p_{1}}({\mathbb{T}}^{2}))_{\theta }=L^{p},$ where $\frac{1}{p}=\frac{1-\theta
 }{p_{0}}+\frac{\theta }{p_{1}}.$
\end{thm}
By duality, because $B.M.O.$ is the dual of $H^{1}$ we get $(L^{2},\
 BMO)_{\theta }=L^{p},$ for $\frac{1}{p}=\frac{1-\theta }{2}.$
 Now by use of a Theorem of Wolff we get:\ \par 

\begin{cor}
\label{dHP2}We have $(L^{1}({\mathbb{T}}^{2}),\ BMO({\mathbb{T}}^{2}))_{\theta
 }=L^{p},$ where $\frac{1}{p}=1-\theta .$
\end{cor}
\quad Hence the assumption (AS) is true for the bidisc.\ \par 

\begin{thm}
\label{iDC8}Let $S$ be a dual bounded sequence in $H^{p}({\mathbb{D}}^{n}),$
 which is also Carleson. Then $S$ is $H^{s}({\mathbb{D}}^{n})$
 interpolating with the LEP, for any $s\in \lbrack 1,p\lbrack
 $ provided  that the assumption (AS) is true.
\end{thm}
\quad In order to prove Theorem~\ref{iDC8}, we shall need some tools.\ \par 

\subsection{BMO functions in the polydisc.}
We shall work in the polydisc ${\mathbb{D}}^{n}$ and we shall
 use the family of kernels, for $p\in \rbrack 1,\infty \lbrack $:\ \par 
\quad \quad \quad $\displaystyle P(a,z):=\frac{(1-\left\vert{a_{1}}\right\vert
 ^{2})^{p-1}}{\left\vert{1-\bar a_{1}z_{1}}\right\vert ^{p}}{\times}\cdot
 \cdot \cdot {\times}\frac{(1-\left\vert{a_{n}}\right\vert ^{2})^{p-1}}{\left\vert{1-\bar
 a_{n}z_{n}}\right\vert ^{p}}.$\ \par 
To see that these kernels are "good" ones, we use the Cayley
 transform ${\mathbb{D}}\rightarrow U:={\mathbb{R}}^{2}_{+}$
 to read them in $U^{n}$:\ \par 
\quad \quad \quad $\displaystyle P(r\alpha ,z)\rightarrow \varphi _{t}(y):=\varphi
 _{t_{1}}(y_{1})\cdot \cdot \cdot \varphi _{t_{n}}(y_{n}):=\frac{1}{t_{1}\cdot
 \cdot \cdot t_{n}}\varphi (\frac{y_{1}}{t_{1}})\cdot \cdot \cdot
 \varphi (\frac{y_{n}}{t_{n}})$\ \par 
with $\varphi (x):=(\frac{1}{x^{2}+1})^{p}.$ Because $\varphi
 $ is ${\mathcal{C}}^{\infty }$ smooth and in $L^{1}({\mathbb{R}})$
 for $p>1,$ $\varphi _{t}(y)$ is a good kernel in the product
 of the upper half planes $U^{n}.$ See ~\cite{ChangFeff85}, ~\cite{Lin86},
 for further details.\ \par 
\quad In particular we shall use the fact that the non tangential maximal
 function $f\rightarrow M_{P}f,$ associated to these kernels
 sends the real $H^{1}({\mathbb{T}}^{n})$ in $L^{1}({\mathbb{T}}^{n})$
 and $L^{p}({\mathbb{T}}^{n})$ in $L^{p}({\mathbb{T}}^{n})$ boundedly.\ \par 
We have the Corollary to the proof of the characterisation by
 A. Chang~\cite{AChang79} of the Carleson measures for the bidisc.
 Let $P_{r}(e^{i\theta },\zeta )=P(re^{i\theta },\zeta )$ be
 a kernel from the family above and set:\ \par 
\quad \quad \quad $\displaystyle u(r_{1}e^{i\theta _{1}},r_{2}e^{i\theta _{2}}):=Pf(re^{i\theta
 }):=P_{r_{2}}\ast (P_{r_{1}}\ast f)(\theta _{1},\theta _{2}).$\ \par 
Let $U\subset {\mathbb{T}}^{2}$ be an open connected set and
 define the region $S(U)$ in ${\mathbb{D}}^{2}$ to be:\ \par 
\quad \quad \quad $\displaystyle S(U):=\lbrace (z_{1},z_{2})\in {\mathbb{D}}^{2}::I_{z_{1}}{\times}I_{z_{2}}\subset
 U\rbrace $\ \par 
where $I_{z_{j}}:=\lbrace e^{i\theta }::\left\vert{\theta -\theta
 _{j,0}}\right\vert <1-r_{j}\rbrace $ with $z_{j}=r_{j}^{i\theta
 _{j,0}}.$ \ \par 
Define as usual for $p\geq 1$ and $f\in L^{p}({\mathbb{T}}^{2}),$\ \par 
\quad \quad \quad $\displaystyle {\mathcal{H}}^{p}({\mathbb{D}}^{2}):=\lbrace Pf::M_{P}f\in
 L^{p}({\mathbb{T}}^{2})\rbrace .$\ \par 

\begin{cor}
Let $\mu $ be a positive measure on ${\mathbb{D}}^{2},$ let $p\geq
 1,$ then $\mu $ is bounded on ${\mathcal{H}}^{p}({\mathbb{D}}^{2})$ iff:\par 
\quad \quad \quad $\mu (S(U))\leq C\left\vert{U}\right\vert $ for all connected,
 open sets $U\subset {\mathbb{T}}^{2}.$
\end{cor}
\quad Proof.\ \par 
As A. Chang~\cite{AChang79} said, it is the same proof as in
 Stein~\cite{SteinDiff70}, page 236. $\hfill\blacksquare $\ \par 
\quad These measures are the Carleson measures of the bidisc.\ \par 
\ \par 
\quad We shall use the following lemma.\ \par 

\begin{lem}
\label{dHP1}Let $P(a,\zeta )$ be a kernel in the above family.
 Let $\omega $ be a Carleson measure in ${\mathbb{D}}^{2}.$ Then
 its balayage by $P(a,\zeta )$ is a $B.M.O.({\mathbb{T}}^{2})$ function.
\end{lem}
\quad Proof.\ \par 
Let\ \par 
\quad \quad \quad $\displaystyle P^{*}\omega (\zeta ):=\int_{{\mathbb{D}}^{2}}{P(a,\zeta
 )d\omega (a)}$\ \par 
be the balayage of $\omega $ by $P.$\ \par 
We know that the dual space of $H^{1}({\mathbb{T}}^{n})$ is $B.M.O.({\mathbb{T}}^{n})$
 by~\cite{DualityH1BMOChang80},\ \par 
~\cite{ChangFeff85}, so it suffices to test $P^{*}\omega (\zeta
 )$ against a smooth function $h$ in ${\mathcal{H}}^{1}({\mathbb{D}}^{2}).$\
 \par 
\quad We have\ \par 
\quad \quad \quad $\displaystyle \int_{{\mathbb{T}}^{2}}{h(\zeta )P^{*}\omega (\zeta
 )d\zeta }=\int_{{\mathbb{T}}^{2}{\times}{\mathbb{D}}^{2}}{h(\zeta
 )P(a,\zeta )d\omega (a)d\zeta }.$\ \par 
Set\ \par 
\quad \quad \quad $\displaystyle H(a):=\int_{{\mathbb{T}}^{2}}{h(\zeta )P(a,\zeta
 )d\zeta },$\ \par 
then we get by the Theorem of Fubini\ \par 
\quad \quad \quad $\displaystyle \int_{{\mathbb{T}}^{2}}{h(\zeta )P^{*}\omega (\zeta
 )d\zeta }=\int_{{\mathbb{D}}^{2}}{H(a)d\omega (a)}.$\ \par 
\quad Then, because $\omega $ is a Carleson measure and $H(a)$ is a
 function in ${\mathcal{H}}^{1}({\mathbb{D}}^{n}),$\ \par 
\quad \quad \quad $\displaystyle \left\vert{\int_{{\mathbb{T}}^{2}}{h(\zeta )P^{*}\omega
 (\zeta )d\zeta }}\right\vert \leq \int_{{\mathbb{D}}^{2}}{\left\vert{H(a)}\right\vert
 d\left\vert{\omega }\right\vert (a)}\leq C{\left\Vert{\omega
 }\right\Vert}_{C}{\left\Vert{h}\right\Vert}_{H^{1}}.$\ \par 
So we get that $P^{*}\omega $ is in $B.M.O.({\mathbb{T}}^{n})$
 with its norm bounded by the Carleson norm of $\omega ,\ C{\left\Vert{\omega
 }\right\Vert}_{C}.$ $\hfill\blacksquare $\ \par 
\ \par 
\quad Now we are in position to prove Theorem~\ref{iDC8}, i.e.\ \par 

\begin{thm}
Let $S$ be a dual bounded sequence in $H^{p}({\mathbb{D}}^{n}),$
 which is also Carleson. Then $S$ is $H^{s}({\mathbb{D}}^{n})$
 interpolating with the LEP, for any $s\in \lbrack 1,p\lbrack
 $ provided  that the assumption (AS) is true.
\end{thm}
\quad Proof.\ \par 
First it is known that all the results above are true for ${\mathbb{D}}^{n},\
 n\geq 2,$ except, up to my knowledge, the result on the interpolation
 between $L^{1}({\mathbb{T}}^{2})$ and $BMO({\mathbb{T}}^{2}).$\ \par 
Suppose that $S$ is a dual bounded sequence in $H^{p}({\mathbb{D}}^{n}),$
 i.e. there is a sequence $\lbrace \rho _{a}\rbrace _{a\in S}\subset
 H^{p}({\mathbb{D}}^{n})$ such that:\ \par 
\quad \quad \quad $\displaystyle \exists C>0,\ \forall a,b\in S,\ {\left\Vert{\rho
 _{a}}\right\Vert}_{p}\leq C,\ {\left\langle{\rho _{a},k_{b,p'}}\right\rangle}=\delta
 _{a,b}.$\ \par 
Now fix $s<p$ and let $\nu \in \ell ^{s}(S).$ We have to show
 that we can interpolate the sequence $\nu $ in $H^{s}({\mathbb{D}}^{n}).$\
 \par 
\quad Let $q$ be such that $\frac{1}{s}=\frac{1}{p}+\frac{1}{q}$ and
 write for $a\in S,\ \nu _{a}=\lambda _{a}\mu _{a}$ with $\lambda
 _{a}:=\frac{\nu _{a}}{\left\vert{\nu _{a}}\right\vert }\left\vert{\nu
 _{a}}\right\vert ^{s/p}$ and $\mu _{a}:=\left\vert{\nu _{a}}\right\vert
 ^{s/q}$; then $\lambda \in \ell ^{p},\ \mu \in \ell ^{q}$ and
 ${\left\Vert{\nu }\right\Vert}_{s}={\left\Vert{\lambda }\right\Vert}_{p}{\left\Vert{\mu
 }\right\Vert}_{q}.$\ \par 
\ \par 
\quad Let, with $\lbrace \gamma _{a}\rbrace _{a\in S}$ to be precised later,\ \par 
\quad \quad \quad $\displaystyle h:=\sum_{a\in S}{\gamma _{a}\nu _{a}\rho _{a}k_{a,q}}=\sum_{a\in
 S}{\gamma _{a}\lambda _{a}\rho _{a}\mu _{a}k_{a,q}}.$\ \par 
Then we get\ \par 
\quad \quad \quad $\displaystyle \forall b\in S,\ h(b)=\sum_{a\in S}{\gamma _{a}\nu
 _{a}\rho _{a}(b)k_{a,q}(b)}=\gamma _{b}\nu _{b}\rho _{b}(b)k_{b,q}(b).$\ \par 
Now we choose $\lbrace \gamma _{a}\rbrace _{a\in S}$ such that
 $\gamma _{b}\nu _{b}\rho _{b}(b)k_{b,q}(b)=\nu _{b}{\left\Vert{k_{b}}\right\Vert}_{s'}.$
 We get:\ \par 
\quad \quad \quad $\displaystyle \gamma _{b}\nu _{b}\rho _{b}(b)k_{b,q}(b)=\gamma
 _{b}{\left\Vert{k_{b}}\right\Vert}_{p'}\frac{{\left\Vert{k_{b}}\right\Vert}_{2}^{2}}{{\left\Vert{k_{b}}\right\Vert}_{q}}=\nu
 _{b}{\left\Vert{k_{b}}\right\Vert}_{s'}$\ \par 
by the choice of $\gamma _{b}.$ Hence the function $h$ interpolates
 the sequence $\nu .$ Moreover $h$ depends linearly on $\nu .$\ \par 
\quad It remains to estimate its $H^{s}({\mathbb{D}}^{n})$ norm.\ \par 
First we estimate $\gamma _{b}.$ We have:\ \par 
\quad \quad \quad $\displaystyle \gamma _{b}\nu _{b}\rho _{b}(b)k_{b,q}(b)=\nu
 _{b}{\left\Vert{k_{b}}\right\Vert}_{s'}\Rightarrow \gamma _{b}=\frac{{\left\Vert{k_{b}}\right\Vert}_{s'}}{\rho
 _{b}(b)k_{b,q}(b)}.$\ \par 
But $1={\left\langle{\rho _{b},\frac{k_{b}}{{\left\Vert{k_{b}}\right\Vert}_{p'}}}\right\rangle}=\frac{\rho
 _{b}(b)}{{\left\Vert{k_{b}}\right\Vert}_{p'}}\Rightarrow \rho
 _{b}(b)={\left\Vert{k_{b}}\right\Vert}_{p'}.$ Also $k_{b,q}(b)=\frac{k_{b}(b)}{{\left\Vert{k_{b}}\right\Vert}_{q}}=\frac{{\left\Vert{k_{b}}\right\Vert}_{2}^{2}}{{\left\Vert{k_{b}}\right\Vert}_{q}}.$
 Hence\ \par 
\quad \quad \quad $\displaystyle \gamma _{b}=\frac{{\left\Vert{k_{b}}\right\Vert}_{s'}{\left\Vert{k_{b}}\right\Vert}_{q}}{{\left\Vert{k_{b}}\right\Vert}_{p'}{\left\Vert{k_{b}}\right\Vert}_{2}^{2}}.$\
 \par 
Now we know that the Structural Hypotheses are true for the polydisc,
 see~\cite{AmarIntInt07} and~\cite{AmarExtInt06}. This means that\ \par 
\quad \quad \quad $\displaystyle \exists \alpha =\alpha _{q}>0,\ \forall a\in {\mathbb{D}}^{n},\
 {\left\Vert{k_{a}}\right\Vert}_{2}^{2}\geq \alpha {\left\Vert{k_{a}}\right\Vert}_{q}{\left\Vert{k_{a}}\right\Vert}_{q'}$\
 \par 
and, for $p,q,s$ such that $1/s=1/p+1/q,$\ \par 
\quad \quad \quad $\displaystyle \exists \beta =\beta _{p,q}>0,\ \forall a\in {\mathbb{D}}^{n},\
 {\left\Vert{k_{a}}\right\Vert}_{s'}\leq \beta {\left\Vert{k_{a}}\right\Vert}_{p'}{\left\Vert{k_{a}}\right\Vert}_{q'}$\
 \par 
Hence, taking $a=b,$ we get $\gamma _{b}\leq \alpha \beta .$\ \par 
Now, we have, using H\"older inequalities\ \par 
\quad \quad \quad $\displaystyle \left\vert{h}\right\vert \leq C\alpha \beta {\left({\sum_{a\in
 S}{\left\vert{\lambda _{a}}\right\vert ^{p}\left\vert{\rho _{a}}\right\vert
 ^{p}}}\right)}^{1/p}{\left({\sum_{a\in S}{\left\vert{\mu _{a}}\right\vert
 ^{p'}\left\vert{k_{a,q}}\right\vert ^{p'}}}\right)}^{1/p'}$\ \par 
with $p'$ the conjugate exponent of $p.$\ \par 
\quad Let $g:={\left({\sum_{a\in S}{\left\vert{\lambda _{a}}\right\vert
 ^{p}\left\vert{\rho _{a}}\right\vert ^{p}}}\right)}^{1/p},$ we have\ \par 
\quad \quad \quad $\displaystyle {\left\Vert{g}\right\Vert}_{p}^{p}=\int_{{\mathbb{T}}^{2}}{\sum_{a\in
 S}{\left\vert{\lambda _{a}}\right\vert ^{p}\left\vert{\rho _{a}}\right\vert
 ^{p}}d\sigma }=\sum_{a\in S}{\left\vert{\lambda _{a}}\right\vert
 ^{p}{\left\Vert{\rho _{a}}\right\Vert}_{p}^{p}}\leq C{\left\Vert{\lambda
 }\right\Vert}_{p}^{p},$\ \par 
hence $g\in L^{p}({\mathbb{T}}^{n}).$ Now we replace $C$ by $C\alpha
 \beta ,$ to ease notation.\ \par 
\quad Let $f:={\left({\sum_{a\in S}{\left\vert{\mu _{a}}\right\vert
 ^{p'}\left\vert{k_{a,q}}\right\vert ^{p'}}}\right)}^{1/p'},$
 we have $\left\vert{h}\right\vert \leq Cfg$ hence\ \par 
\quad \quad \quad $\displaystyle {\left\Vert{h}\right\Vert}_{s}^{s}=\int_{{\mathbb{T}}^{n}}{\left\vert{h}\right\vert
 ^{s}d\sigma }\leq C^{s}\int_{{\mathbb{T}}^{n}}{f^{s}g^{s}d\sigma }.$\ \par 
Applying H\"older inequalities with exponents $p/s,\ q/s$ because
 $\frac{1}{s}=\frac{1}{p}+\frac{1}{q},$ we get\ \par 
\quad \quad \quad $\displaystyle {\left\Vert{h}\right\Vert}_{s}^{s}\leq C^{s}{\left\Vert{g}\right\Vert}_{p}^{s}{\left\Vert{f}\right\Vert}_{q}^{s}.$\
 \par 
To have $h\in L^{s}({\mathbb{T}}^{n}),$ it remains to prove that\ \par 
\quad \quad \quad $\displaystyle f:={\left({\sum_{a\in S}{\left\vert{\mu _{a}}\right\vert
 ^{p'}\left\vert{k_{a,q}}\right\vert ^{p'}}}\right)}^{1/p'}\in
 L^{q}({\mathbb{T}}^{n}).$\ \par 
So let $F:=f^{p'}=\sum_{a\in S}{\left\vert{\mu _{a}}\right\vert
 ^{p'}\left\vert{k_{a,q}}\right\vert ^{p'}},$ we shall show that
 $F\in L^{q/p'}({\mathbb{T}}^{n}),\ (q>p'\Rightarrow q/p'>1).$\ \par 
\quad From\ \par 
\quad \quad \quad $\displaystyle k_{a,q}(z)=\frac{(1-\left\vert{a_{1}}\right\vert
 ^{2})^{1/q'}}{(1-\bar a_{1}z_{1})}{\times}\cdot \cdot \cdot
 {\times}\frac{(1-\left\vert{a_{n}}\right\vert ^{2})^{1/q'}}{(1-\bar
 a_{n}z_{n})}.$\ \par 
we get\ \par 
\quad \quad \quad $\displaystyle \left\vert{k_{a,q}}\right\vert ^{p'}=\frac{(1-\left\vert{a_{1}}\right\vert
 ^{2})^{p'/q'}}{\left\vert{1-\bar a_{1}z_{1}}\right\vert ^{p'}}{\times}\cdot
 \cdot \cdot {\times}\frac{(1-\left\vert{a_{n}}\right\vert ^{2})^{p'/q'}}{\left\vert{1-\bar
 a_{n}z_{n}}\right\vert ^{p'}}=$\ \par 
\quad \quad \quad \quad \quad \quad \quad $\displaystyle =(1-\left\vert{a_{1}}\right\vert ^{2})^{-p'/q}\cdot
 \cdot \cdot (1-\left\vert{a_{n}}\right\vert ^{n})^{-p'/q}{\times}$\ \par 
\quad \quad \quad \quad \quad \quad \quad \quad \quad $\displaystyle {\times}\frac{(1-\left\vert{a_{1}}\right\vert
 ^{2})^{p'-1}}{\left\vert{1-\bar a_{1}z_{1}}\right\vert ^{p'}}{\times}\cdot
 \cdot \cdot {\times}\frac{(1-\left\vert{a_{n}}\right\vert ^{2})^{p'-1}}{\left\vert{1-\bar
 a_{n}z_{n}}\right\vert ^{p'}}{\times}(1-\left\vert{a_{1}}\right\vert
 ^{2})\cdot \cdot \cdot (1-\left\vert{a_{n}}\right\vert ^{n}).$\ \par 
\quad We make the hypothesis that the measure\ \par 
\quad \quad \quad $\displaystyle d\omega (z):=\sum_{a\in S}{(1-\left\vert{a_{1}}\right\vert
 ^{2})\cdot \cdot \cdot (1-\left\vert{a_{n}}\right\vert ^{2})\delta
 _{a}}$\ \par 
is a Carleson measure in ${\mathbb{D}}^{n}.$ Then we shall apply
 the Lemma~\ref{dHP1} to get that the balayage of $d\omega (z)$ by\ \par 
\quad \quad \quad $\displaystyle P(a,z):=\frac{(1-\left\vert{a_{1}}\right\vert
 ^{2})^{p'-1}}{\left\vert{1-\bar a_{1}z_{1}}\right\vert ^{p'}}{\times}\cdot
 \cdot \cdot {\times}\frac{(1-\left\vert{a_{n}}\right\vert ^{2})^{p'-1}}{\left\vert{1-\bar
 a_{n}z_{n}}\right\vert ^{p'}}$\ \par 
is in $B.M.O.({\mathbb{T}}^{n}),$ because this kernel verifies
 the conditions of Lemma~\ref{dHP1}.\ \par 
\quad Now if $\mu $ is a bounded measure in ${\mathbb{D}}^{n}$ then
 its balayage by $P(a,z)$ is in $L^{1}({\mathbb{T}}^{n})$ because
 $P(a,z)$ is of norm $1$ in $L^{1}({\mathbb{T}}^{n})$  for any
 $a$ fixed in ${\mathbb{D}}^{n}.$\ \par 
\ \par 
\quad Hence if $\gamma (z)\in L^{t}(d\omega )$ by interpolation using
 the assumption (AS), valid for the bidisc by Lin's Theorem~\ref{AI10},
 we get that\ \par 
\quad \quad \quad $\displaystyle P^{*}(\gamma d\omega )(\zeta )\in L^{t}({\mathbb{T}}^{n}).$\
 \par 
\quad So it remains to see that\ \par 
\quad \quad \quad $\displaystyle \gamma :=\left\vert{\mu _{a}}\right\vert ^{p'}(1-\left\vert{a_{1}}\right\vert
 ^{2})^{-p'/q}\cdot \cdot \cdot (1-\left\vert{a_{n}}\right\vert
 ^{2})^{-p'/q}\in L^{q/p'}(d\omega ).$\ \par 
But\ \par 
\quad \quad \quad $\displaystyle \int_{{\mathbb{D}}^{n}}{\gamma ^{q/p'}d\omega
 }=\sum_{a\in S}{\left\vert{\mu _{a}}\right\vert ^{q}(1-\left\vert{a_{1}}\right\vert
 ^{2})^{-1}\cdot \cdot \cdot (1-\left\vert{a_{n}}\right\vert
 ^{n})^{-1}{\times}(1-\left\vert{a_{1}}\right\vert ^{2})\cdot
 \cdot \cdot (1-\left\vert{a_{n}}\right\vert ^{n})},$\ \par 
hence\ \par 
\quad \quad \quad $\displaystyle \int_{{\mathbb{D}}^{n}}{\gamma ^{q/p'}d\omega
 }=\sum_{a\in S}{\left\vert{\mu _{a}}\right\vert ^{q}}={\left\Vert{\mu
 }\right\Vert}_{\ell (S)}^{q}<\infty .$\ \par 
So we prove that $F\in L^{q/p'}({\mathbb{T}}^{n}),$ hence $f:={\left({\sum_{a\in
 S}{\left\vert{\mu _{a}}\right\vert ^{p'}\left\vert{k_{a,q}}\right\vert
 ^{p'}}}\right)}^{1/p'}\in L^{q}({\mathbb{T}}^{n})$ and we are
 done. $\hfill\blacksquare $\ \par 

\section{Appendix. Assumption (AS)}
\quad We shall see that the assumption (AS) would be a corollary of
 a weak factorisation.\ \par 
\quad We shall suppose that\ \par 
\quad \quad \quad $\exists C>0,\ \forall 1\leq p\leq 2,\ \forall f\in H^{p}({\mathbb{D}}^{n}),\
 \exists g_{j},h_{j}\in H^{2p}({\mathbb{D}}^{n})::f=\sum_{j\in
 {\mathbb{N}}}{g_{j}h_{j}}$\ \par 
and\ \par 
\quad \quad \quad $\sum_{j\in {\mathbb{N}}}{{\left\Vert{g_{j}}\right\Vert}_{2p}{\left\Vert{h_{j}}\right\Vert}_{2p}}\leq
 C{\left\Vert{f}\right\Vert}_{p},$\ \par 
i.e. we shall suppose that $H^{p}({\mathbb{D}}^{n})$ has \emph{the
 extended weak factorisation property} (EWF) for $1\leq p\leq
 2$ with a uniform constant $C.$\ \par 
This true for $p=1$ by M. T. Lacey and E. Terwilleger~\cite{LaceyTerWi09}
  and for $p>1$ but the constant is a priori not uniform.\ \par 
\ \par 
Because the Szeg\"o projection $P:L^{2}\rightarrow H^{2}$ is
 uniformly bounded on $L^{p}({\mathbb{T}}^{n}),$ for $2\leq p\leq
 4$ we get:\ \par 

\begin{lem}
\label{AI11}We have:\par 
\quad \quad \quad $\displaystyle (H^{2}({\mathbb{D}}^{n}),\ H^{4}({\mathbb{D}}^{n}))_{\theta
 }=H^{p}({\mathbb{D}}^{n})$ with $\displaystyle \frac{1}{p}=\frac{1-\theta
 }{2}+\frac{\theta }{4}.$
\end{lem}
\quad Proof.\ \par 
First take $f\in H^{p}({\mathbb{D}}^{n}).$ Then $f\in L^{p}({\mathbb{T}}^{n})$
 and we have $Pf=f.$ We already know that $(L^{2}({\mathbb{T}}^{n}),L^{4}({\mathbb{T}}^{n}))_{\theta
 }=L^{p}({\mathbb{T}}^{n})$ so it exists a holomorphic function
 $F(z,\zeta )$ for $z$ in the strip $0<\Re z<1$ and in $L^{2}({\mathbb{T}}^{n})$
 in $\zeta $ on $\Re z=0$ and in $L^{4}({\mathbb{T}}^{n})$ in
 $\zeta $ on $\Re z=1$   such that:\ \par 
\quad \quad \quad $\displaystyle F(\theta ,\cdot )=f,\ \exists C>0,\ \forall y\in
 {\mathbb{R}},\ {\left\Vert{F(iy,\cdot )}\right\Vert}_{L^{2}}\leq
 C{\left\Vert{f}\right\Vert}_{L^{p}},\ \ {\left\Vert{F(1+iy,\cdot
 )}\right\Vert}_{L^{4}}\leq C{\left\Vert{f}\right\Vert}_{L^{p}}.$\ \par 
Set $G(z,\zeta ):=P_{\zeta }F(z,\cdot )(\zeta ),$ then $G(z,\zeta
 )$ is holomorphic as a function of $\zeta $ and, because $2\leq
 p\leq 4,$ there is a uniform $C>0$ such that:\ \par 
\quad \quad \quad $\displaystyle G(\theta ,\zeta )=P_{\zeta }F(\theta )(\zeta )=Pf(\zeta
 )=f(\zeta ),\ \ \forall y\in {\mathbb{R}},\ {\left\Vert{G(iy,\cdot
 )}\right\Vert}_{L^{2}}\leq C{\left\Vert{f}\right\Vert}_{L^{p}},\
 \ {\left\Vert{G(1+iy,\cdot )}\right\Vert}_{L^{4}}\leq C{\left\Vert{f}\right\Vert}_{L^{p}}.$\
 \par 
Because $P$ is linear, we keep the fact that $G(z,\zeta )$ is
 also holomorphic in $z$ in the strip $0<\Re z<1.$ So we get
 that $f\in (H^{2}({\mathbb{D}}^{n}),H^{4}({\mathbb{D}}^{n}))_{\theta
 }$ hence $(H^{2}({\mathbb{D}}^{n}),H^{4}({\mathbb{D}}^{n}))_{\theta
 }\supset H^{p}({\mathbb{D}}^{n}).$\ \par 
\quad To prove the converse, we take $f\in (H^{2}({\mathbb{D}}^{n}),H^{4}({\mathbb{D}}^{n}))_{\theta
 }\subset (L^{2}({\mathbb{T}}^{n}),L^{4}({\mathbb{T}}^{n}))_{\theta
 }=L^{p}({\mathbb{T}}^{n})$ hence it remains to see that $f$
 is holomorphic. We use the fact that:\ \par 
\quad \quad \quad $\displaystyle f(e^{i\psi })=F(\theta ,e^{i\psi })=\int_{{\mathbb{R}}}{P_{\theta
 }(iy)F(iy,e^{i\psi })dy}+\int_{{\mathbb{R}}}{P_{\theta }(1+iy)F(1+iy,e^{i\psi
 })dy}$\ \par 
where $P_{z}(\cdot )$ is the Poisson kernel of the strip $0<\Re
 z<1.$ To see that $f$ extends holomorphically in ${\mathbb{D}}^{n},$
 we compute its Fourier coefficients, with $\psi =(\psi _{1},...,\psi
 _{n}),\ e^{i\psi }:=e^{i\psi _{1}}{\times}\cdot \cdot \cdot
 {\times}e^{i\psi _{n}},\ d\psi :=d\psi _{1}\cdot \cdot \cdot
 d\psi _{n}$:\ \par 
\quad $\displaystyle \hat f(k):=\int_{{\mathbb{T}}^{n}}{f(e^{i\psi
 })e^{ik\psi }d\psi }$\ \par 
\quad \quad \quad $\displaystyle =\int_{{\mathbb{T}}^{n}}{{\left({\int_{{\mathbb{R}}}{P_{\theta
 }(iy)F(iy,e^{i\psi })dy}+\int_{{\mathbb{R}}}{P_{\theta }(1+iy)F(1+iy,e^{i\psi
 })dy}}\right)}e^{ik\psi }d\psi }.$\ \par 
Using the Theorem of Fubini we get\ \par 
\quad \quad \quad $\displaystyle \hat f(k)=\int_{{\mathbb{R}}}{P_{\theta }(iy)\hat
 F(iy,k)dy}+\int_{{\mathbb{R}}}{P_{\theta }(1+iy)\hat F(1+iy,k)dy}$\ \par 
hence, if $k=(k_{1},...,k_{n})$ contains a negative entry, then
 $\hat F(iy,k)=\hat F(1+iy,k)=0$ because these are holomorphic
 functions in $\zeta .$ So we get that, if $k=(k_{1},...,k_{n})$
 contains a negative entry, then $\hat f(k)=0$ hence $f$ extends
 holomorphically in ${\mathbb{D}}^{n}.$ The proof is complete.
 $\hfill\blacksquare $\ \par 
\ \par 
\quad Then we have\ \par 

\begin{thm}
\label{iC1}Suppose that $H^{p}({\mathbb{D}}^{n})$ has the extended
 weak factorisation property for $1\leq p\leq 2$ with a uniform
 constant $C.$ Then we have by the complex method: $(H^{1},H^{2})_{\theta
 }=H^{p}$ with $\frac{1}{p}=1-\frac{\theta }{2}.$
\end{thm}
\quad Proof.\ \par 
Take $f\in (H^{1},H^{2})_{\theta }$ with $\frac{1}{p}=1-\frac{\theta
 }{2}$ then we get easily that $f\in H^{p}.$ The point is to
 prove that if $f\in H^{p}$ then $f\in (H^{1},H^{2})_{\theta
 }$ with $\frac{1}{p}=1-\frac{\theta }{2}.$\ \par 
\quad So let $f\in H^{p}.$ By assumption $f=\sum_{j\in {\mathbb{N}}}{g_{j}h_{j}}$
 with $g_{j},h_{j}\in H^{2p}.$ Because $p\geq 1,$ we get $2p\geq
 2$ hence, by Lemma~\ref{AI11}, there is a holomorphic function
 $G_{j}(z)$ in the strip $0<\Re z<1$ such that $G_{j}(\theta
 )=g_{j}$ and ${\left\Vert{G_{j}(iy)}\right\Vert}_{H^{2}}\leq
 C{\left\Vert{g_{j}}\right\Vert}_{H^{2p}}$ and ${\left\Vert{G_{j}(1+iy)}\right\Vert}_{H^{4}}\leq
 C{\left\Vert{g_{j}}\right\Vert}_{H^{2p}}.$ The same there is
 a holomorphic function $H_{j}(z)$ in the strip $0<\Re z<1$ such
 that $H_{j}(\theta )=h_{j}$ and ${\left\Vert{H_{j}(iy)}\right\Vert}_{H^{2}}\leq
 C{\left\Vert{g_{j}}\right\Vert}_{H^{2p}}$ and ${\left\Vert{H_{j}(1+iy)}\right\Vert}_{H^{4}}\leq
 C{\left\Vert{h_{j}}\right\Vert}_{H^{2p}}.$\ \par 
\quad Let $F(z):=\sum_{j\in {\mathbb{N}}}{G_{j}(z)H_{j}(z)}.$ Then
 $F(z)$ is a holomorphic function in the strip $0<\Re z<1$ and\ \par 
\quad \quad \quad $\displaystyle {\left\Vert{F(iy)}\right\Vert}_{H^{1}}\leq \sum_{j\in
 {\mathbb{N}}}{{\left\Vert{G_{j}(iy)}\right\Vert}_{H^{2}}{\left\Vert{H_{j}(iy)}\right\Vert}_{H^{2}}}\leq
 C{\left\Vert{f}\right\Vert}_{p},$\ \par 
and\ \par 
\quad \quad \quad $\displaystyle {\left\Vert{F(1+iy)}\right\Vert}_{H^{2}}\leq \sum_{j\in
 {\mathbb{N}}}{{\left\Vert{G_{j}(1+iy)}\right\Vert}_{H^{4}}{\left\Vert{H_{j}(1+iy)}\right\Vert}_{H^{4}}}\leq
 C{\left\Vert{f}\right\Vert}_{p}.$\ \par 
Moreover we have\ \par 
\quad \quad \quad $\displaystyle F(\theta )=\sum_{j\in {\mathbb{N}}}{G_{j}(\theta
 )H_{j}(\theta )}=\sum_{j\in {\mathbb{N}}}{g_{j}h_{j}}=f.$\ \par 
Hence we proved that $f\in (H^{1},H^{2})_{\theta }$ with $\frac{1}{p}=1-\frac{\theta
 }{2}.$ The proof is complete. $\hfill\blacksquare $\ \par 
\ \par 
\quad Now we set $BMO({\mathbb{T}}^{n})$ the \emph{real} $BMO$ space
 and $H^{1}_{R}$ the \emph{real} $H^{1}$ space. Then\ \par 

\begin{cor}
We have $(L^{1},\ BMO)_{\theta }=L^{p},$ where $\frac{1}{p}=1-\theta .$
\end{cor}
\quad Proof.\ \par 
First Theorem~\ref{iC1} gives\ \par 
\quad \quad \quad $\displaystyle (H_{R}^{1},L^{2})_{\theta }=L^{p}$ with $\displaystyle
 \frac{1}{p}=1-\frac{\theta }{2}$\ \par 
because if $f\in L^{p},\ p>1,$ then $f\in H_{R}^{p}$ hence from
 Theorem~\ref{iC1} we get that $f\in (H_{R}^{1},H_{R}^{2})_{\theta
 }$ so, a fortiori, $f\in (H_{R}^{1},L^{2})_{\theta }.$ If $f\in
 (H_{R}^{1},L^{2})_{\theta },$ then $f\in (L^{1},L^{2})_{\theta
 }=L^{p}$ and we are done.\ \par 
\quad By duality, because $B.M.O.$ is the dual of $H_{R}^{1}$ by~\cite{AChang79},
 we get\ \par 
\quad \quad \quad \quad $\displaystyle (L^{2},BMO)_{\theta }=L^{p}$ with $\displaystyle
 \frac{1}{p}=\frac{1-\theta }{2}.$\ \par 
Now choose $q>2.$ We shall use the extrapolation Theorem by T.
 Wolff~\cite{Wolff81} which says that if:\ \par 
\quad \quad \quad $\displaystyle (L^{1},L^{q})_{\theta }=L^{r}$ and $\displaystyle
 (L^{2},BMO)_{\theta }=L^{s}$\ \par 
then\ \par 
\quad \quad \quad $\displaystyle (L^{1},\ BMO)_{\theta }=L^{p},$ where $\displaystyle
 \frac{1}{p}=1-\theta .$\ \par 
The proof is complete. $\hfill\blacksquare $\ \par 
\ \par 
\quad We know that the weak factorisation is true for $H^{p}({\mathbb{D}}^{n}),\
 p>1$ and for $H^{1}({\mathbb{D}}^{n}).$ An "heuristic" proof
 of its validity for $H^{p}({\mathbb{D}}^{n}),\ p\geq 1,$ with
 a uniform constant, is the following.\ \par 
\quad We take the proof for $H^{1}({\mathbb{D}}^{n})$ done by M. T.
 Lacey and E. Terwilleger~\cite{LaceyTerWi09} which consists
 on a careful study of a function $b$ in $BMO$ and we replace
 $BMO$ by $H^{p'}({\mathbb{D}}^{n})$ the dual space of $H^{p}({\mathbb{D}}^{n}).$
 Then we keep careful track of the constants and, hopefully,
 we get the weak factorisation with uniform constants between
 $1\leq p\leq 2.$ Of course this is just heuristic!\ \par 
\ \par 

\bibliographystyle{/usr/local/texlive/2017/texmf-dist/bibtex/bst/base/apalike}

\end{document}